\documentclass[12pt]{article}
\usepackage{amssymb,amsthm,bm}

\title{Abelian functions satisfy an Algebraic Addition Theorem}

\author{Mark B. Villarino\\
        Escuela de Matem\'atica, Universidad de Costa Rica,\\
        2060 San Jos\'e, Costa Rica}

\date{April 16, 1998}

\theoremstyle{plain}
\newtheorem{thm}{Theorem}

\theoremstyle{definition}
\newtheorem{defn}{Definition}

\newcommand{\f}{\varphi}              

\newcommand{\C}{\mathbb{C}}           

\newcommand{\uu}{\bm{u}}              
\newcommand{\vv}{\bm{v}}              

\newcommand{\x}{\times}               

\newcommand{\row}[3]{{#1}_{#2},\dots,{#1}_{#3}} 

\begin{document}

\maketitle

\begin{abstract}
We prove that $n$ independent abelian functions admit an algebraic 
addition theorem, with no appeal to theta functions.
\end{abstract}


\vspace{1pc}

A famous and classical theorem of Weierstrass states:

\begin{thm}[{\rm Weierstrass, 1869}]
\label{th:abelianAAT}
Any $n$ analytically independent Abelian functions in $n$ complex 
variables that belong to the same $2n$-dimensional period lattice 
satisfy an algebraic addition theorem.
\end{thm}

All published proofs of this theorem depend on the 
Weierstrass--Riemann ``Theta-satz'' which affirms that any abelian 
function is globally representable as the quotient of two 
(generalized) theta functions~\cite{SiegelL}.

We show in this note that the theorem has nothing to do with theta 
functions \textit{per~se}, but is a simple (indeed, almost trivial) 
consequence of another, even more famous, theorem of 
Weierstrass~\cite{Weierstrass}, Thimm~\cite{Thimm} and
Siegel~\cite{SiegelM} on algebraic dependence.

\begin{thm}[{\rm Weierstrass--Thimm--Siegel}]
\label{th:WThS}
Any $n + 1$ meromorphic functions on a compact complex space of
dimension~$n$ are algebraically dependent.
\end{thm}

\begin{defn}
\label{df:AAT}
A set of $n$ meromorphic functions of $n$ complex variables, 
$\f_k(\row{u}{1}{n})$, $k = 1,\dots,n$, satisfy an 
\textbf{algebraic addition theorem} if and only if there exist $n$ 
polynomials $\row{G}{1}{n}$ in $2n + 1$ complex variables with 
complex coefficients, such that the equations
$$
G_k\bigl[ \f_k(\uu + \vv); \f_1(\uu),\dots,\f_n(\uu); 
          \f_1(\vv),\dots,\f_n(\vv) \bigr] = 0
$$
for $k = 1,\dots,n$, hold identically for all 
$\uu = (\row{u}{1}{n}) \in \C^n$ and $\vv = (\row{v}{1}{n}) \in \C^n$ 
where the $\f_i$ are defined.
\end{defn}

\begin{defn}
The meromorphic functions $\row{\f}{1}{k}$ on the compact complex 
space $X$ are called \textbf{algebraically dependent} iff there 
exists a polynomial $Q(\row{z}{1}{k}) \not\equiv 0$ with complex 
coefficients such that 
$$
Q\bigl( \f_1(\uu),\dots,\f_k(\uu) \bigr) = 0
$$
holds identically for all $\uu \in X$ where the $\f_k$ are defined.
\end{defn}

\begin{proof}[Proof of Theorem~\ref{th:abelianAAT}]
Let $X$ be the period parallelotope for the lattice $\Omega$ of 
periods of the abelian functions $\f_1(\uu),\dots,\f_n(\uu)$. Then 
the cartesian product $X \x X$ is a compact complex connected 
manifold of dimension~$2n$.

The function $\f_k(\uu + \vv)$ is a meromorphic function of 
$2n$~variables $(\row{u}{1}{n},\allowbreak\row{v}{1}{n})$ on $X \x X$,
as are the $2n$ functions $\f_k(\uu)$ and $\f_k(\vv)$. Therefore, the 
$2n + 1$ meromorphic functions $\f_k(\uu + \vv)$; 
$\f_1(\uu),\dots,\f_n(\uu)$; $\f_1(\vv),\dots,\f_n(\vv)$ on the 
$2n$-dimensional compact complex space~$X \x X$ are necessarily 
algebraically dependent, by Theorem~\ref{th:WThS}. This conclusion 
holds for each $\f_k(\uu + \vv)$, $k = 1,\dots,n$. But, by 
Definition~\ref{df:AAT}, this means that $\row{\f}{1}{k}$ satisfy an 
algebraic addition theorem.
\end{proof}


\end{document}